\documentclass[10pt]{article}
\usepackage{amsmath,amssymb,amsthm}
\usepackage{geometry}
\newtheorem{Theorem}{Theorem}[section]
\newtheorem{Definition}{Definition}[section]
\newtheorem{Lemma}[Definition]{Lemma}

\numberwithin{equation}{section}
\allowdisplaybreaks
\usepackage[colorlinks=true,linkcolor=blue,anchorcolor=blue,citecolor=blue]{hyperref}

\begin{document}

\pagenumbering{arabic} \setcounter{page}{1}
\begin{center}
{\large \bf Large-time behavior of solutions to the 1D  nematic liquid crystal with $\gamma >1$}\\
${\rm Xiuqing\; Wang^{1*}}\quad {\rm Man \; Liu^{1}}   \quad{\rm Yuming\; Qin^{2}}$\\
 1. Institute for Differential Equation Theory and Applications, Faculty of Science,  \\Kunming University of Science and
Technology,
Kunming, P. R. China.\\
2. School of Mathematics and Statistics,  Department of Mathematics,  Institute for Nonlinear Sciences,\\ Donghua University, Shanghai,  P. R. China.\\

E-mail:   daqingwang@kust.edu.cn \quad  liuman@stu.kust.edu.cn   \quad  yuming@dhu.edu.cn\\

\end{center}

\begin{abstract}~~
We study the large-time behavior of  solutions to the system modeling the compressible flow of nematic liquid crystals with a polytropic equation of state $P(\rho)=\rho^{\gamma}$ $(\gamma>1)$ in Lagrangian coordinates. In contrast to the isothermal case $\gamma=1$ treated in \cite{4,5},  we establish a new representation formula for the specific volume (Lemma~\ref{lem:2.3}) and, using Jensen's inequality, derive time-uniform estimates that yield global upper and lower bounds for the specific volume. A key ingredient is the reconstruction of the $D$--$Z$ structure for $\gamma>1$, which gives an exponential decay of $Z(t)$ and thus overcomes the difficulty caused by the power-type pressure.

\end{abstract}

\section{Introduction}
The nematic liquid crystal flow \cite{1,2,6} is described by the following coupled system:
\begin{subequations}\label{1.1}
\begin{align}
&\rho_{t}+(\rho u)_{x}=0, \\[3pt]
&(\rho u)_{t}+(\rho u^{2})_{x}+\bigl(P(\rho)\bigr)_{x}
   = u_{xx}- (|n_{x}|^{2})_{x}, \\[3pt]
&n_{t}+u\,n_{x}=  n_{xx}+|n_{x}|^{2}n .
\end{align}
\end{subequations}
Mathematically, the one-dimensional compressible nematic liquid crystal system exhibits a hybrid structure: the mass equation is hyperbolic, the momentum equation is parabolic-hyperbolic, and the orientation equation is a linear parabolic equation with an algebraic constraint $|n|=1$.
 Here $\rho(x, t) $ denotes the density,
$u$ and $n$ are unknown
functions standing for the   velocity  field,  macroscopic average of the  orientation field, respectively.

We consider the  initial-boundary data for  system (\ref{1.1}) in the domain $\{ (x,t) | 0<x<1, t\in [0,+\infty)\}$
\begin{equation}\label{1.2}
 \begin{cases}
u(0,t)=u(1,t)=0, n_{x}(0,t)=n_{x}(1,t)=0,\\
(\rho, u,n )(\cdot, 0)=(\rho_{0},u_{0},n_{0}),
 \end{cases}
\end{equation}
and the  polytropic gas law $ P(\rho)= \rho^{\gamma}(\gamma > 1). $
Under the Lagrangian coordinates, i.e., $ y=\int_{0}^{x}\rho(\xi,t)d\xi, $
system (\ref{1.1})-(\ref{1.2}) transformed  into the following:
\begin{subequations}\label{1.3}
\begin{align}
&v_{t}=u_{y},  \\[3pt]
&u_{t}=\Bigl(-P+\frac{u_{y}}{v}-\frac{|n_{y}|^{2}}{v^{2}}\Bigr)_{y},
 \\[3pt]
&n_{t}=\frac{1}{v}\Bigl(\frac{n_{y}}{v}\Bigr)_{y}+\frac{|n_{y}|^{2}}{v^{2}}\,n,
\end{align}
\end{subequations}
\begin{equation}\label{1.4}
 \begin{cases}
u(0,t)=u(1,t)=0, n_{y}(0,t)=n_{y}(1,t)=0,\\
(v, u,n )(\cdot, 0)=(v_{0},u_{0},n_{0}),
 \end{cases}
\end{equation}
where $v=1/\rho$ and $v_{0}=1/\rho_{0}$.

In recent years, there have been numerous papers concerning the existence theory of solutions to the liquid crystal system, see References \cite{7,12,13,14}.
For the one-dimensional case, Ding-Wang-Wen \cite{3} established in 2011 the global existence of weak solutions in $H^{1}\times H^{1}_{0}\times H^{2}$ under the isothermal pressure law $P(\rho)=\rho$.
Qin and Huang \cite{4} established the global high-regularity theory for system (\ref{1.1})-(\ref{1.2}) and obtained an exact representation formula for the specific volume $v$, a key tool for analyzing large-time behavior. Subsequently. Qin-Feng \cite{5} derived time-uniform estimates and the large-time asymptotics in the isothermal case $\gamma = 1$.
However, the isothermal pressure law is physically restrictive, corresponding only to the isothermal process of an ideal gas in a heat bath. The more physically relevant pressure law is the polytropic gas law  which describes the adiabatic process of real gases. The extension of Qin-Feng \cite{5} to the polytropic regime $(P(\rho)=\rho^\gamma$ with $\gamma>1)$  remains an open and non-trivial problem.

This paper extends these results to polytropic gases,  $(P(\rho)=\rho^\gamma$ with $\gamma>1)$.
Our main contribution is to carry the large-time analysis over to $\gamma>1$: by combining the representation formula of Lemma \ref{lem:2.3} with energy estimates and Jensen's inequality, we obtain time-uniform estimates that yield global upper and lower bounds for $v$ and the associated large-time control. The key ingredient is a $\gamma$-dependent reconstruction of the $v, D,Z$ representation which yields exponential decay of $Z(t)$, thereby resolving the nonlinear difficulties caused by the power-type pressure.

\begin{Theorem}\label{thm:1}
Suppose the initial data $(v_0, u_0, n_0)\in H^{i}(0,1)\times H^{i}_{0}(0,1)\times H^{i+1}(0,1) (i=1,2)$ satisfying  the compatibility conditions. Then the system (\ref{1.3})-(\ref{1.4}) possesses a unique global smooth solution $(v (t), u (t), n (t))\in H^{i}(0,1)\times H^{i}_{0}(0,1)\times H^{i+1}(0,1)$,  and there exists a constant $C>0$, such that
\begin{equation}\label{1.5}
0<C_{1}^{-1}\leq v(y,t)\leq C ,\quad\forall\,(y,t)\in[0,1]\times[0,+\infty),
\end{equation}
and
 \begin{equation} \label{1.6}
\|v(t) \|_{H^{i}}^{2}+ \|u(t)\|_{H^{i}}^{2}+\|n(t) \|_{H^{i+1}}^{2}+\int_{0}^{\infty}\left(\|u_{y}(t) \|_{H^{i }}^{2}+\|n_{y}(t) \|_{H^{i+1}}^{2}\right)(\tau)d\tau \leq C .
 \end{equation}
Moreover, as $t\rightarrow +\infty$, we have
\begin{equation}\label{1.7}
\|v(t)-\bar{v}\|_{H^{i}}\rightarrow 0,\quad \|u(t)\|_{H^{i}}\rightarrow 0,\quad \|n(t)-\bar{n}\|_{H^{i+1}}\rightarrow 0,
\end{equation}
where $\bar{v}=\int_0^1v(y,t)\,dy=\bar v_0, \bar{n}=\lim\limits_{t\to\infty}\int_0^1n(y,t)\,dy$.
\end{Theorem}

\section{Uniform estimates}

\begin{Lemma}\label{lem:2.1}

Under the assumptions in Theorem \ref{thm:1}, the following estimates
hold in Lagrangian coordinates, for any $t>0$,
\begin{align}
&\int _{0}^{1}\Bigl[u^{2}+\frac{v^{1-\gamma}}{\gamma -1}+\frac{|n_{y}|^{2}}{v}\Bigr]
   (y,t)\,\mathrm{d}y+\int _{0}^{t}\int_{0}^{1}\Bigl(\frac{u_{y}^{2}}{v}
   +\frac{1}{v}\Bigl|n\times\left(\frac{n_{y}}{v}\right)_{y}\Bigr|^{2}\Bigr)
   \mathrm{d}y\mathrm{d}s\leq C_{1}, \label{eq:lem3-1}\\[3pt]
&\int_{0}^{1}\frac{|n_{y}|^{2}}{v}\,\mathrm{d}y
   +\int_{0}^{t}\int_{0}^{1} \frac{1}{v}\Bigl[\Bigl(\frac{n_{y}}{v}\Bigr)_{y}\Bigr]^{2}
   \mathrm{d}y\mathrm{d}s\leq C_{1}. \label{eq:lem3-2}
\end{align}
\end{Lemma}

\begin{Lemma}\label{lem:2.3}
Let \begin{align}
U(y,t)=\int_0^y u(\xi,t)d\xi, \label{2.5}
\qquad
\bar v_0=\int_0^1v_0(y)dy,
\end{align}
and define
\begin{align}
D(y,t)
&=
v_0(y)
\exp\left\{
U(y,t)-U_0(y)
-\frac1{\bar v_0}
\left[
\int_0^1vUdy-\int_0^1v_0U_0dy
\right]
\right\},
  \label{2.6}\\[3pt]
Z(t)
&=
\exp\left\{
-\frac1{\bar v_0}
\int_0^t\int_0^1
\left(
u^2+v^{1-\gamma}
+\frac{|n_y|^2}{v}
\right)dyds
\right\}.
  \label{2.7}
\end{align}
Then  for any $t\geq 0$,  the solution $v(y,t)$ to problems (\ref{1.3})-(\ref{1.4})
possesses the following expression:
\begin{align}
v^\gamma(y,t)
={}D^\gamma(y,t)Z^\gamma(t)
\exp\left\{
\gamma\int_0^t
\frac{|n_y(y,s)|^2}{v^2(y,s)}ds
\right\}\times
\left[
1+
\gamma\int_0^t
D^{-\gamma}(y,s)Z^{-\gamma}(s)
\exp\left\{
-\gamma\int_0^s
\frac{|n_y(y,\tau)|^2}{v^2(y,\tau)}d\tau
\right\}ds
\right],
\end{align}
and
\begin{align}
C^{-1}\le D(y,t)\le C.
\end{align}
\end{Lemma}

\begin{proof}
Since $
u(0,t)=u(1,t)=0,$ integrate equation (\ref{1.3}a) over \(y\),
there holds
\begin{equation}
\int_0^1v(y,t)\,dy=\bar v_0.
\label{2.10}
\end{equation}
(1) We first prove that \(D\) is uniformly bounded. By the Cauchy-Schwarz inequality and Lemma \ref{lem:2.1}, we get
\begin{equation}
\left|
\int_0^1vU\,dy
\right|
\le
\|U\|_{L^\infty}\int_0^1v\,dy\le
\|u(t)\|_{L^2}\int_0^1v\,dy
\le C\bar v_0,\quad \text{similarly,} \left|
\int_0^1v_0U_0\,dy
\right|
\le C.
\label{2.12}
\end{equation}
Thus there exists a constant \(C_D>1\), such that
\begin{equation}
C_D^{-1}\le D(y,t)\le C_D.
\label{2.13}
\end{equation}
(2) We next prove the exponential decay rate of \(Z(t)\).
 Define
\begin{equation}
G(t)=
\int_0^1
\left(
u^2+v^{1-\gamma}
+\frac{|n_y|^2}{v}
\right)dy,
\label{2.14}
\end{equation}
by Lemma \ref{lem:2.1}, we have
\begin{equation}
G(t)\le C,
\label{2.15}\quad \text{and}\quad Z(t)\ge e^{-C_Zt}.
\end{equation}

Note that $ G(t)\ge \int_0^1v^{1-\gamma}dy.\ $
Since \(\gamma>1\), the function $ f(v)=v^{1-\gamma} $
satisfies $ f''(v)=\gamma(\gamma-1)v^{-\gamma-1}>0 $,\(f\) is convex.
Using Jensen's inequality and mass conservation equation (\ref{2.10}),
we obtain
\begin{equation}
G(t)\ge\int_0^1v^{1-\gamma}dy
\ge
\left(\int_0^1v\,dy\right)^{1-\gamma}\ge \bar v_0^{\,1-\gamma}
=:c_0>0.
\label{2.17}
\end{equation}
According to (\ref{2.7}), for any \(0\le s\le t\), it holds
\begin{equation}
e^{-C_Z(t-s)}
\le
\frac{Z(t)}{Z(s)}=
\exp\left[
-\frac1{\bar v_0}\int_s^tG(\tau)d\tau
\right]
\le
e^{-c_Z(t-s)}, \quad  C_Z=\frac{C}{\bar v_0},c_Z=\frac{c_0}{\bar v_0}.
\label{2.18}
\end{equation}

(3) Derive the representation formula for \(v\).  Let
\begin{equation}
A=
\frac{u_y}{v}
-v^{-\gamma}
-\frac{|n_y|^2}{v^2},
\label{2.19}
\end{equation}
then the momentum equation  (\ref{1.3}b)  reads
\begin{equation}
u_t=A_y.
\label{2.20}
\end{equation}
According to (\ref{2.5}) and (\ref{2.20}),
we have
\begin{equation}
U_t
=
\int_0^yu_t(\xi,t)d\xi
=
A(y,t)-A(0,t).
\label{2.21}
\end{equation}

On the one  hand, using equation (\ref{1.3}a) and (\ref{2.19}), we get
\begin{equation}
(\ln v)_t
=\frac{v_t}{v}
=
A+
v^{-\gamma}
+\frac{|n_y|^2}{v^2}.
\label{2.22}
\end{equation}
Multiplying (\ref{2.21}) by \(v\), and integrate over \([0,1]\),
\begin{equation}
\int_0^1vA\,dy
=
A(0,t)\bar v_0
+
\int_0^1vU_t\,dy.
\label{2.23}
\end{equation}
For the left-hand side of (\ref{2.23}), by equation (\ref{2.19}), we have
\begin{equation}
\begin{aligned}
\int_0^1vA\,dy
=
\int_0^1
\left(
u_y
-v^{1-\gamma}
-\frac{|n_y|^2}{v}
\right)dy
=
-\int_0^1v^{1-\gamma}dy
-\int_0^1\frac{|n_y|^2}{v}dy.
\end{aligned}
\label{2.24}
\end{equation}
Next we consider the right-hand side of (\ref{2.23}), since
\[
\frac{d}{dt}\int_0^1vUdy
=
\int_0^1v_tUdy+\int_0^1vU_tdy=
\int_0^1u_yUdy+\int_0^1vU_tdy,
\]
using integrate by parts and $U_y=u$,
we obtain $ \int_0^1u_yUdy = -\int_0^1u^2dy,\ $
then
\begin{equation}
\int_0^1vU_tdy
=
\frac{d}{dt}\int_0^1vUdy
+
\int_0^1u^2dy.
\label{2.25}
\end{equation}
Substitute \eqref{2.23}--\eqref{2.25}:
\[
A(0,t)\bar v_0
+
\frac{d}{dt}\int_0^1vUdy
+
\int_0^1u^2dy
=
-\int_0^1v^{1-\gamma}dy
-\int_0^1\frac{|n_y|^2}{v}dy,
\]
therefore
\begin{equation}
A(0,t)
=
-\frac1{\bar v_0}
\frac{d}{dt}\int_0^1vUdy
-\frac1{\bar v_0}G(t).
\label{2.26}
\end{equation}

By the definition of \(D, Z\), and equation (\ref{2.21}),
it follows that
\begin{equation}
\frac{D_t}{D}
=
U_t
-
\frac1{\bar v_0}
\frac{d}{dt}\int_0^1vUdy
=
A-A(0,t)
-
\frac1{\bar v_0}
\frac{d}{dt}\int_0^1vUdy, \quad \frac{Z_t}{Z}
=
-\frac{G(t)}{\bar v_0},
\label{2.27}
\end{equation}
and added equation (\ref{2.22}), using \eqref{2.26}, so we get
\begin{equation}
\frac{d}{dt}
\ln\frac{v}{DZ}=A+v^{-\gamma}+\frac{|n_y|^2}{v^2}
-\left[
A-A(0,t)
-\frac1{\bar v_0}\frac{d}{dt}\int vU
\right]
+\frac{G}{\bar v_0}
=
v^{-\gamma}
+\frac{|n_y|^2}{v^2}.
\label{2.29}
\end{equation}

On the other hand,
set
\[
W(y,t)=\frac{v(y,t)}{D(y,t)Z(t)},
\]
then
\begin{equation}
W_t
=
\left(
v^{-\gamma}
+\frac{|n_y|^2}{v^2}
\right)W.
\label{2.30}
\end{equation}
and we have $ W(y,0)=1. $
Let $ X=W^\gamma, $
By using \eqref{2.30}, since
$
v^{-\gamma}X
=
(DZ)^{-\gamma}, $
therefore
\begin{equation}
X_t
=
\gamma(DZ)^{-\gamma}
+
\gamma\frac{|n_y|^2}{v^2}X,
\label{2.31}
\end{equation}
this is a linear ODE of $X$!
Therefore
\begin{equation}
\begin{aligned}
W^\gamma(y,t)
={}
\exp\left\{
\gamma\int_0^t
\frac{|n_y(y,s)|^2}{v(y,s)^2}ds
\right\}
\times
\left[
1+
\gamma\int_0^t
(D(y,s)Z(s))^{-\gamma}
\exp\left\{
-\gamma\int_0^s
\frac{|n_y(y,\tau)|^2}{v(y,\tau)^2}d\tau
\right\}ds
\right].
\end{aligned}
\label{2.32}
\end{equation}
This yields the Lemma \ref{lem:2.3}:
\[
\begin{aligned}
v^\gamma(y,t)
={}
D^\gamma(y,t)Z^\gamma(t)
e^{\gamma\int_0^t\frac{|n_y|^2}{v^2}ds}
\times
\left[
1+
\gamma\int_0^t
D^{-\gamma}(y,s)Z^{-\gamma}(s)
e^{-\gamma\int_0^s\frac{|n_y|^2}{v^2}d\tau}
\,ds
\right].
\end{aligned}
\]
\end{proof}

\begin{Lemma}\label{lem:2.4}
There holds
\begin{equation}\label{eq:lem5-1}
0<C_{1}^{-1}\leq v(y,t)\leq C_{1},\quad\forall\,(y,t)\in[0,1]\times[0,+\infty),
\end{equation}
\begin{equation}\label{eq:lem5-2}
\int_{0} ^{t}\bigl\| n_{y}(s) \bigr\|^{2}_{L^{\infty}}\,\mathrm{d}s\leq C_{1},
\quad\forall\,t>0.
\end{equation}
\end{Lemma}

\begin{proof}
(1) First note that $ e^{\gamma\int_0^t\frac{|n_y|^2}{v^2}ds}\ge1. $
Thus by \eqref{2.32},we have
$ W^\gamma(y,t)
\ge
\gamma\int_0^t(DZ)^{-\gamma}(y,s)\,ds, $
hence
\[
v^\gamma(y,t)
\ge
\gamma D^\gamma(y,t)Z^\gamma(t)
\int_0^tD^{-\gamma}(y,s)Z^{-\gamma}(s)\,ds.
\]
Take the last unit time interval. If \(t\ge1\), using (\ref{2.13}) and (\ref{2.18}) , we obtain
\[
\left(\frac{D(t)}{D(s)}\right)^\gamma
\ge C_D^{-2\gamma} ,\quad v^\gamma(y,t)
\ge
\gamma
\int_{t-1}^t
\left(
\frac{D(y,t)}{D(y,s)}
\right)^\gamma
\left(
\frac{Z(t)}{Z(s)}
\right)^\gamma ds\ge
\gamma C_D^{-2\gamma}e^{-\gamma C_Z},
\]
thus
\begin{equation}
v(y,t)\ge c_*>0,
\qquad t\ge1.
\label{2.35}
\end{equation}
For \(0\le t\le1\), using \eqref{2.30}, we have $ W(y,t)\ge1. $
Thus
\[
v(y,t)=D(y,t)Z(t)W(y,t)
\ge C_D^{-1}e^{-C_Zt}
\ge C_D^{-1}e^{-C_Z}.
\]
Combining the above estimates, we obtain
\begin{equation}
v(y,t)\ge C_1^{-1},
\qquad
(y,t)\in[0,1]\times[0,\infty).
\label{2.36}
\end{equation}

(2) Let $ f=\frac{n_y}{v}. $
By the boundary conditions (\ref{1.4}),
we have $ f(0,t)=f(1,t)=0. $
According to the  Sobolev inequality yields
\begin{equation}
\|f(t)\|_{L^\infty}^2
\le
C\left(
\|f(t)\|_{L^2}^2
+
\|f_y(t)\|_{L^2}^2
\right)=C\left( \int_0^1\frac{|n_y|^2}{v^2}dy+ \int_0^1
\left|
\left(\frac{n_y}{v}\right)_y
\right|^2dy
\right).
\label{2.37}
\end{equation}
Since \(v\ge C_1^{-1}\), and by Lemma \ref{lem:2.1}, we have
\begin{equation}
\int_0^t
\left\|
\frac{n_y}{v}
\right\|_{L^\infty}^2ds
\le C_1
\qquad\forall t>0.
\label{2.39}
\end{equation}

(3) Using  \eqref{2.39} and \eqref{2.32}, we obtain
\[
\begin{aligned}
v^\gamma(y,t)
\le {}&
D^\gamma(y,t)Z^\gamma(t)e^{\gamma C_1} +
\gamma D^\gamma(y,t)Z^\gamma(t)
\int_0^t
D^{-\gamma}(y,s)Z^{-\gamma}(s)
e^{\gamma C_1}\,ds.
\end{aligned}
\]
Using the upper and lower bounds of \(D\), and  the upper bound in \eqref{2.18},
then we lead to
\[
v^\gamma(y,t)
\le
C+
C\int_0^te^{-\gamma c_Z(t-s)}ds.
\]
Note that
$
\int_0^te^{-\gamma c_Z(t-s)}ds
\le
\frac1{\gamma c_Z}, $
therefore
$
v^\gamma(y,t)\le C, $
namely
\begin{equation}
v(y,t)\le C.
\label{2.40}
\end{equation}
Combined with \eqref{2.36},  \eqref{2.39}, the proof is complete.
\end{proof}

\section{Higher-order estimates and long-time behaviors}
In this section, based on the previous basic energy estimates and uniform bounds of \(v\) we further establish the uniform higher-order estimates of the solution and prove its long-time asymptotic behavior.   Lemma \ref{lem:2.1}-\ref{lem:2.4} serve as the foundation for the subsequent higher-order estimates.

\begin{Lemma}\label{lem:3.1}
 Under the assumptions of Theorem \ref{thm:1}.
Let
$
(v,u,n)
$
be a smooth solution to the system (\ref{1.3})-(\ref{1.4}),  there exists a constant \(C>0\), such that
\begin{align}
\sup_{t\geq0}
\left(
\|v(t)\|_{H^1}^2
+\|u(t)\|_{H^1}^2
+\|n(t)\|_{H^{2}}^2
\right)+\int_0^\infty
\left(
\|u_y(t)\|_{H^1}^2
+\|n_y(t)\|_{H^{2}}^2
\right)\,dt
\leq C.
\end{align}

\end{Lemma}

\begin{proof}
Using (\ref{1.3})a, we rewrite (\ref{1.3})b as
\begin{equation}
\left(u-\frac{v_{y}}{v}\right)_t=\frac{\gamma v_{y}}{v^{\gamma +1}}+\frac{2|n_{y}|^{2} v_{y}}{v^{3}}-\frac{2 n_{y}\cdot n_{yy}}{v^{2}},
\label{3.2}
\end{equation}
multiplying the equation  (\ref{3.2}) by $\left(u-\frac{v_{y}}{v}\right)$, integrate over $(0,1)\times [0,t]$, and using Lemma \ref{lem:2.1}-\ref{lem:2.4}, we deduce
\begin{align}
  \int_{0}^{1}\left|u-\frac{v_{y}}{v}\right|^{2}(t) dy
+\int_{0}^{t}\int_{0}^{1}\left(v_{y}^{2}+\frac{|n_{y}|^{2} v_{y}^{2}}{v^{4}}+\frac{|n_{yy}|^{2} }{v^{3}}\right) dy ds
\leq C_{1}.
\label{3.3}
\end{align}

Multiplying the equation  (\ref{1.3})b  by $u_{yy}$, integrate over $(0,1)$ about $y$, we deduce
\begin{align}
\frac{\mathrm{1}}{\mathrm{2}}\frac{d}{dt} \int_{0}^{1}|u_{y}|^{2}dy +\int\nolimits^{\mathrm{1}}_{\mathrm{0}}{\frac{\left|u_{yy}\right|^{\mathrm{2}}}{v}} dy=-\gamma  \int\nolimits^{\mathrm{1}}_{\mathrm{0}}{\frac{v_y u_{yy}}{v^{\gamma+\mathrm{1}}}} dy+\int\nolimits^{\mathrm{1}}_{\mathrm{0}}{\frac{u_y v_y u_{yy}}{v^{\mathrm{2}}}} dy+\mathrm{2}\int\nolimits^{\mathrm{1}}_{\mathrm{0}}{\frac{(n_y\cdot n_{yy})u_{yy}}{v^{\mathrm{2}}}} dy
-2\int\nolimits^{\mathrm{1}}_{\mathrm{0}}{\frac{|n_y|^{\mathrm{2}} v_y u_{yy}}{v^{\mathrm{3}}}} dy,
\label{3.4}
\end{align}
according to $\|u\|_{L^{\infty}}\leq C(\|u\|_{L^{2}}+\|u_y\|_{L^{2}})$ and $a^{\frac{3}{2}}b\leq \varepsilon a^{2}+C_{\varepsilon} b^{4}$,  integrate (\ref{3.4}) over $(0,t)$ about $s$,  we obtain
 \begin{align}
& \int_{0}^{1}|u_{y}|^{2}dy +\int_{0}^{t}\int\nolimits^{\mathrm{1}}_{\mathrm{0}}{\frac{\left|u_{yy}\right|^{\mathrm{2}}}{v}} dyds \leq C.
\label{3.5}
\end{align}

Multiplying the equation  (\ref{1.3})c  by $v^{2}$,
taking $\|\cdot\|_{L^{2}}$, we get
\begin{equation}
\|n_{yy}\|_{L^{2}}\leq C_{1}\left(\|n_{t}\|_{L^{2}}+1\right) . \label{3.7}
\end{equation}

Take the $t$  derivative of the equation (\ref{1.3})c,
\begin{equation}
n_{tt}=\frac{\mathrm{1}}{v}\left(\frac{n_{ty}}{v}\right)_y- \frac{u_y}{v^{\mathrm{2}}}\left(\frac{n_y}{v}\right)_y+\frac{\mathrm{2}(n_y\cdot n_{ty})}{v^{\mathrm{2}}} n-\frac{\mathrm{2}|n_y|^{\mathrm{2}} u_y}{v^{\mathrm{3}}} n+\frac{\left|n_y\right|^{\mathrm{2}}}{v^{\mathrm{2}}} n_t,
\label{3.8}
\end{equation}
 multiply the  equations by  $n_{t}$,  we have
\begin{equation}
 \int_{0}^{1}|n_t(t)|^{2}dy +\int_{0}^{t}\int_{0}^{1} | n_{ty}(t)|^{2}dyds \leq C_{1}.
\label{3.9}
\end{equation}

Take the $y$  derivative of the equation (\ref{1.3})c,
\begin{equation}
n_{yyy}=\mathrm{2} v v_y n_t+v^{\mathrm{2}} n_{ty}+\frac{n_{yy} v_y+n_y u_{yy}}{v}- \frac{n_y v^{\mathrm{2}}_y}{v^{\mathrm{2}}}-\mathrm{2}(n_y\cdot n_{yy})n- |n_y|^{\mathrm{2}} n_y,
\label{3.10}
\end{equation}
taking $\|\cdot\|_{L^{2}}$, and using the above inequalities, we get
\begin{equation}
\int_{0}^{t}\|n_{yyy}\|_{L^{2}}^{2}ds \leq C_{1} . \label{3.11}
\end{equation}
This completes the proof of the lemma.

\end{proof}

\begin{Lemma}\label{lem:3.2}
 Under the assumptions of Theorem \ref{thm:1}.
Let
$
(v,u,n)
$
be a smooth solution to the system (\ref{1.3})-(\ref{1.4}), there exists a constant \(C>0\), such that
\begin{align}
\sup_{t\geq0}
\left(
\|v(t)\|_{H^2}^2
+\|u(t)\|_{H^2}^2
+\|n(t)\|_{H^{3}}^2
\right)+\int_0^\infty
\left(
\|u_y(t)\|_{H^2}^2
+\|n_y(t)\|_{H^{3}}^2
\right)\,dt
\leq C,
\end{align}
and \begin{align}
\lim_{t\to\infty}\|v(t)-\bar{v}\|_{H^i}=0,
\qquad
\lim_{t\to\infty}\|u(t)\|_{H^i}=0,
\qquad
\lim_{t\to\infty}\|n(t)-\bar{n}\|_{H^{i+1}}=0,
\qquad i=1,2.
\label{3.25}
\end{align}
\end{Lemma}

\begin{proof}
By virtue of the boundary conditions, using Lemmas 8-11 in \cite{5} and Lemmas 12-21 in \cite{5},   we obtain the assertion. This completes the proof of the lemma.
\end{proof}

\begingroup
\small
\setlength{\baselineskip}{11pt}

\endgroup

\end{document}